\documentclass[11pt]{article}

\usepackage{amssymb,amsmath}

\usepackage{eepic}
\usepackage[dvipdfm]{pict2e}

\usepackage{color}
\usepackage[dvipdfm]{graphicx}

\newtheorem{theorem}{Theorem}[section]

\newtheorem{Lemma}[theorem]{Lemma}
\newtheorem{prop}[theorem]{Proposition}

\numberwithin{equation}{section}

\newcommand{\non}{\nonumber}
\newcommand{\cR}{{\mathbb R}}

\newcommand{\eq}[1]{\mbox{\rm {(\ref{#1})}}}

\title{\Large\bf Traveling waves for models of phase transitions
of  solids  driven by configurational forces  }

\author{\small
\sc Shuichi Kawashima$^1$\thanks{E-mail:
kawashim@math.kyushu-u.ac.jp}  \ \ and\ \ Peicheng
Zhu$^{2,3}$\thanks{E-mail: zhu@bcamath.org}
\\ \\ \small
 $^1$Graduate School of Mathematics,  Kyushu University \\
\small  Fukuoka 812-8581,\ \ Japan\\
\small  and\\
\small  $^2$Basque Center for Applied Mathematics\\
\small Building 500, Bizkaia Technology Park \\
\small E-48160 Derio,\ \ Spain\\
\small  $^3$IKERBASQUE, Basque Foundation for Science\\
\small  E-48011 Bilbao,\ \ Spain}

\date{ }

\begin{document}

\maketitle

\centerline{\bf Abstract}
\vskip0.2cm This article is concerned with
the existence of traveling wave solutions, including standing waves,
to some  models based on configurational forces, describing
respectively the diffusionless phase transformations of solid
materials, e.g., Steel, and phase transitions due to interface
motion by interface diffusion, e.g., Sintering. These models were
proposed by Alber and Zhu in \cite{Alber04}. We consider both the
order-parameter-conserved case and the non-conserved one, under
suitable assumptions. Also we compare our results with the
corresponding ones  for the Allen-Cahn and the Cahn-Hilliard
equations coupled with linear elasticity, which are models for
diffusion-dominated phase transformations in elastic solids.


\section{Introduction}

In his work \cite{Alber00}, H.-D. Alber started the study of models
of sharp interface describing phase transitions in alloys. Such
models are of Hamiltonian type and   complex from both theoretical
and numerical viewpoints. Thus in \cite{Alber04,Alber06}, as an
approximation of those models, the phase field models were derived
from the corresponding sharp interface ones for describing phase
transitions in elastically deformable solids. The material under
consideration possesses microstructures, and its phases are
characterized by the structure of the crystal lattice, in which the
atoms are arranged. An interface between different material phases
moves if the crystal lattice is transformed from one structure to
the other across   the interface.

Often phase transformations are triggered by diffusion processes.
Two well known models for diffusion dominated transformations are
the Allen-Cahn and Cahn-Hilliard equations
 (see \cite{Allen79,Cahn59,Cahn59a}).

In this article  we shall consider two models, one describing
diffusionless transformations (also called martensitic
transformations, cf. \cite[p. 162]{Hornbogen}) of solid materials,
e.g., Steel, and another  for phase transitions due to interface
motion by interface diffusion, e.g., Sintering.  The derivation of
these two models can be found in \cite{Alber04,Alber06}. Both models
are based on a model for phase transitions driven by configurational
forces, see \cite{AK90}. For more details on configurational forces,
we refer the reader to the book by Gurtin \cite{Gurtin00}, and the
references cited therein. The models studied here are also related
to the work by Mullins \cite{Mullins} on the development of grooves.
Some mathematical results have been obtained for special cases of
those two types of models, see \cite{Alber04,Alber06,Alber08} for
the existence of weak solutions, and \cite{Alber08a} for the
justification of the diffusive models by using asymptotic expansion
methods.

In the present paper, we study the    existence of traveling wave
solutions, including standing wave solutions, for the phase field
models  in \cite{Alber04,Alber06}. The models are formulated for
materials with two different phases, say, phase $1$ and phase $2$,
in which the order parameter $v=v(x,t)\in \cR$ is assumed to take
values near $v_-$ or $v_+$ if the material is in the phase $1$ or
$2$, respectively, where we assume, without loss of generality, that
 $v_-<v_+$, and $x=(x_1,\cdots,x_n)\in \cR^n$ denotes the material
point and $t\geq 0$ is the time variable. Let $u=(u^i(x,t))\in
\cR^n$ be the displacement of the material. Then the first model in
which the order parameter is not conserved is formulated as
\begin{eqnarray}
&-\,{\rm div}_x\,T=b,
\label{I1a}   \\[0.2cm]
&v_t=-\,(\psi_v(\varepsilon,v)-\Delta_x v)|\nabla_x v|.
\label{I2a}
\end{eqnarray}
The second model is   as follows
\begin{eqnarray}
&-\,{\rm div}_x\,T=b,
  \label{II1a}    \\[0.2cm]
&v_t = {\rm div}_x\big\{\nabla_x (\psi_v(\varepsilon,v)
-\Delta_x v)|\nabla_x v|\big\}.
  \label{II2a}
\end{eqnarray}
In the latter,   the order parameter is conserved. Eq. \eq{II2a} can
be written in the form of a conservation law:
$$
v_t = {\rm div}_x \,  {\bf J},
$$
where $  {\bf J}$ is  a flux defined by
$$  {\bf J} = \nabla_x (\psi_v(\varepsilon,v) -\Delta_x
v)|\nabla_x v|.
$$
%
However \eq{I2a} can not be rewritten in the conservative form.
%

We give now a more detailed description of the above systems. The
symbols ${\rm div}_x$, $\nabla_x$ and $\Delta_x$ denote the
divergence, the gradient and the Laplacian with respect to
$x\in\cR^n$. The $\varepsilon=(\varepsilon_{ij})\in{\mathcal S}^n$
denotes the strain tensor which is given by
\begin{equation}
\varepsilon_{ij}=\frac{1}{2}(u^i_{x_j}+u^j_{x_i}),
\end{equation}
where ${\mathcal S}^n$ denotes the set of all $n\times n$ real
symmetric matrices. $T=(T_{ij})\in{\mathcal S}^n$ is the Cauchy
stress tensor and is given explicitly in terms of $(\varepsilon,v)$,
see \eq{CR1} below. Also, $\psi(\varepsilon,v)$ is a part of the
free energy $\psi^*=\psi(\varepsilon,v)+\frac{1}{2}|\nabla_xv|^2$,
which is given in \eq{CR2} below, and  $\psi_v(\varepsilon,v)$
denotes the derivative of $\psi(\varepsilon,v)$ with respect to $v$.
Finally, $b=b(x,t)$ denotes the volume force which is assumed to be
a given function of $(x,t)$.

To close the systems, we need to add the following constitutive
relations:

\bigskip
\noindent{\bf (CR1)}
\begin{equation}
 T=D(\varepsilon-\bar\varepsilon(v)) ,
 \label{CR1}
\end{equation}

\bigskip
\noindent{\bf (CR2)}
\begin{equation}
 \psi(\varepsilon,v)= \frac{1}{2}\,   D(\varepsilon-\bar\varepsilon (v))
  :(\varepsilon-\bar\varepsilon(v)) + f(v)  .
 \label{CR2}
\end{equation}

\bigskip
\noindent
Here $D=(D_{kl}^{ij})$ is the elasticity tensor which is assumed to
be a linear symmetric positive definite operator from ${\mathcal S}^n$
into ${\mathcal S}^n$; this means that for any
$\sigma=(\sigma_{ij})\in {\mathcal S}^n$,
\begin{equation*}
D\sigma
  =\big(\sum_{k,l=1}^nD^{ij}_{kl}\sigma_{kl}\big)
  \in{\mathcal S}^n,
\qquad
(D\sigma):\sigma\geq c\,|\sigma|^2
\end{equation*}
for some positive constant $c$. Here we denote the scalar product
of two matrices $\sigma$ and $\tau$ by
$\sigma:\tau=\sum_{i,j=1}^n\sigma_{ij}\tau_{ij}$, and the
corresponding norm by $|\sigma|=(\sigma:\sigma)^{1/2}$.
Also, $\bar\varepsilon(v)\in{\mathcal S}^n$ is assumed to be linear
in $v$:
\begin{equation}\label{1.6}
\bar\varepsilon(v)=\varepsilon^{(0)}+\varepsilon^{(1)}v,
\end{equation}
where $\varepsilon^{(0)}=(\varepsilon^{(0)}_{ij})\in{\mathcal S}^n$
and $\varepsilon^{(1)}=( \varepsilon^{(1)}_{ij})\in {\mathcal S}^n$
are transformation tensors which are assumed to be constant
matrices, and $f(v)$ is a given smooth double-well potential with
two local minima at $v = v_-$ and $v = v_+$. Consequently, the pair
$(v,u)(x,t)$ or $(v,\varepsilon)(x,t)$ is the unknown of our systems
\eqref{I1a} -- \eqref{I2a} and \eqref{II1a} -- \eqref{II2a}.

We can slightly rewrite our systems. A simple computation, using
\eq{CR1} and \eq{CR2}, yields
\begin{equation*}
\psi_v(\varepsilon,v)
=-\,D(\varepsilon-\bar\varepsilon(v)):\varepsilon^{(1)}+f'(v)
=-\,T:\varepsilon^{(1)}+f'(v),
\end{equation*}
where $f'$ denotes the derivative of $f$.
Therefore the systems become
\begin{eqnarray}
&-\,{\rm div}_x\,T= b ,
  \label{I1}   \\[0.2cm]
&v_t=\big(T:\varepsilon^{(1)}-f'(v)+\Delta_xv\big)|\nabla_x\,v|,
 \label{I2}
\end{eqnarray}
and
\begin{eqnarray}
&-\,{\rm div}_x\,T=b ,
  \label{II1}  \\[0.2cm]
&v_t=-\,{\rm div}_x\big\{\nabla_x
  \big(T:\varepsilon^{(1)}-f'(v)+\Delta_xv\big)|\nabla_xv|\big\},
  \label{II2}
\end{eqnarray}
respectively. The evolution equations \eq{I2} and \eq{II2} for the
order parameter $v$ are non-uniformly parabolic because of the
degenerate terms $|\nabla_xv|\Delta_xv$ and ${\rm
div}_x(\nabla_x\Delta_xv|\nabla_xv|)$, respectively.

\smallskip
The main purpose of this article is to study the existence of
traveling wave   and standing wave solutions to the systems \eq{I1}
-- \eq{I2} and \eq{II1}  --  \eq{II2}. To this end, we assume that
the volume force $b=b(x,t)$ is identically zero, i.e., $b\equiv 0$,
and look for one-dimensional solutions of the form
\begin{equation*}
v=v(x_1,t), \qquad u=(u^1(x_1,t),0,\cdots,0).
\end{equation*}
Under suitable structural assumptions, it will be observed that
our $n$-dimensional systems are essentially reduced to the
one-dimensional equations
\begin{equation}
\label{1.11}
v_t=(\mu-f'(v)+v_{x_1x_1})|v_{x_1}|
\end{equation}
and
\begin{equation}\label{1.12}
v_t=-\{(-f'(v)+v_{x_1x_1})_{x_1}|v_{x_1}|\}_{x_1},
\end{equation}
respectively, where $\mu$ is a constant to be determined.
 Moreover, we will show that these one-dimensional
equations \eqref{1.11} and \eqref{1.12} have a traveling wave
solution $v=v(x_1- s t)$ and a standing wave solution $v=v(x_1)$,
respectively, both connecting the phase states $v_-$ and $v_+$,
provided that $f(v_+)=f(v_-)$. A comparison in Section~4, will show
that this existence result of standing wave solutions to \eq{1.12}
is the same to the known results for the Cahn-Hilliard equation,
however if $\mu\not=0$, results of traveling waves to  \eq{1.11} are
quite different from those for Allen-Cahn equation.

The remaining parts of this article are as follows: In Section~2,
under suitable structural assumptions, we shall reduce the
multi-dimensional systems into the one-dimensional ones, which can
be finally reduced to the one single equation \eqref{1.11} or
\eqref{1.12} for each case. Then we shall investigate, in Section~3,
the existence of traveling wave  and standing wave solutions for
these reduced one-dimensional systems. The last section is devoted
to the interesting comparison of our results with the existing ones
for the Allen-Cahn or Cahn-Hilliard equation coupled with linear
elasticity.


\section{Reduction to one-dimensional equations}
\setcounter{equation}{0}

In this section, under suitable structural assumptions, we reduce
our multi-dimensional systems \eqref{I1}  --  \eqref{I2} and
\eqref{II1}  --  \eqref{II2} into the corresponding ones in the
one-dimensional form. To this end, we first assume that the volume
force $b=b(x,t)$ is identically zero, namely,
\begin{equation}
b\equiv 0,
\end{equation}
and look for one-dimensional solutions of the form
\begin{equation}\label{2.2}
v=v(x_1,t), \qquad u=(u^1(x_1,t),0,\cdots,0).
\end{equation}
We now calculate the Cauchy stress tensor $T=(T_{ij})$ for these
particular solutions. For $u=(u^1(x_1,t),0,\cdots,0)$, the strain
tensor $\varepsilon=(\varepsilon_{ij})$ has the simplest form
\begin{eqnarray}
\varepsilon_{11}=u^1_{x_1}, \qquad
\varepsilon_{ij}=0 \ \ {\rm for}\ \ (i,j)\neq(1,1).
\end{eqnarray}
Consequently, we obtain from \eq{CR1} and \eqref{1.6} that
\begin{equation}\label{2.4}
\begin{split}
T_{ij}
&=D^{ij}:\big(\varepsilon
  -(\varepsilon^{(0)}+\varepsilon^{(1)}v)\big) \\[1mm]
&=D^{ij}_{11}u^1_{x_1}-(D^{ij}:\varepsilon^{(0)})
  -(D^{ij}:\varepsilon^{(1)})v,
\end{split}
\end{equation}
where $D^{ij}=(D^{ij}_{kl})\in{\mathcal S}^n$.
In particular, we have
\begin{equation}\label{2.5}
T_{11}=D^{11}_{11}u^1_{x_1}
  -(D^{11}:\varepsilon^{(0)})-(D^{11}:\varepsilon^{(1)})v.
\end{equation}
Eliminating $u^1_{x_1}$ in \eq{2.4} by using \eq{2.5} yields an
expression of $T$ in terms of $T_{11}$ and $v$.
Namely, we obtain
\begin{equation}\label{2.6}
T=\sigma T_{11}+\tau^{(0)}+\tau^{(1)}v,
\end{equation}
where $\sigma=(\sigma_{ij})$, $\tau^{(0)}=(\tau^{(0)}_{ij})$ and
$\tau^{(1)}=(\tau^{(1)}_{ij})$ are the constant tensors in
${\mathcal S}^n$ given explicitly as
\begin{equation}\label{2.7}
\begin{split}
\sigma_{ij}
&=D_{11}^{ij}/D_{11}^{11}, \\[1mm]
\tau^{(m)}_{ij}
&=(1/D_{11}^{11})\big\{D_{11}^{ij}(D^{11}:\varepsilon^{(m)})
  -D_{11}^{11}(D^{ij}:\varepsilon^{(m)})\big\},\quad m=0,1.
\end{split}
\end{equation}
In fact, we have from \eqref{2.4} and \eqref{2.5} that
\begin{equation*}
\begin{split}
T_{ij}
&=(D^{ij}_{11}/D^{11}_{11})\big\{T_{11}
  +(D^{11}:\varepsilon^{(0)})+( D^{11}:\varepsilon^{(1)})v\big\}
  -(D^{ij}:\varepsilon^{(0)})-( D^{ij}:\varepsilon^{(1)})v \\[1mm]
&=\sigma_{ij}T_{11}+\tau_{ij}^{(0)}+\tau_{ij}^{(1)}v,
\end{split}
\end{equation*}
which gives the expression \eqref{2.6}.
Notice that $\sigma_{11}=1$ and $\tau^{(0)}_{11}=\tau^{(1)}_{11}=0$.
Also, using \eqref{2.6}, we see that
\begin{equation}\label{2.8}
T:\varepsilon^{(1)}=\alpha T_{11}+\beta+\gamma v,
\end{equation}
where $\alpha$, $\beta$ and $\gamma $ are the constants given by
\begin{equation}\label{2.9}
\begin{split}
&\alpha=\sigma:\varepsilon^{(1)}
  =(1/D_{11}^{11})(D_{11}:\varepsilon^{(1)}), \\[1mm]
&\beta=\tau^{(0)}:\varepsilon^{(1)}
  =(1/D_{11}^{11})\big\{
  (D_{11}:\varepsilon^{(1)})(D^{11}:\varepsilon^{(0)})
  -D_{11}^{11}(D\varepsilon^{(0)}:\varepsilon^{(1)})\big\}, \\[1mm]
&\gamma=\tau^{(1)}:\varepsilon^{(1)}
  =(1/D_{11}^{11})\big\{
  (D_{11}:\varepsilon^{(1)})(D^{11}:\varepsilon^{(1)})
  -D_{11}^{11}(D\varepsilon^{(1)}:\varepsilon^{(1)})\big\}
\end{split}
\end{equation}
with $D_{11}=(D^{ij}_{11})\in{\mathcal S}^n$.
These observations are summarized as follows.

\begin{Lemma}\label{Lem2.1}
Let $u=(u^1(x_1,t),0,\cdots,0)$. Then the Cauchy stress tensor
$T=(T_{ij})$ has the form \eqref{2.6}, in which the component
$T_{11}$ is given by \eqref{2.5}. Moreover, we have \eqref{2.8}.
\end{Lemma}

Now, we assume the following structural conditions:

\bigskip
\noindent {\bf (A1)}
\begin{eqnarray}
D^{11}_{11}(D^{i1}:\varepsilon^{(1)})
=D^{i1}_{11}(D^{11}:\varepsilon^{(1)}) \ \ {\rm for}\ \
 i=1,\cdots,n ,
 \label{A1}
\end{eqnarray}

\bigskip
\noindent {\bf (A2)}
\begin{eqnarray}
(D_{11}:\varepsilon^{(1)})(D^{11}:\varepsilon^{(1)})
=D^{11}_{11}(D\varepsilon^{(1)}:\varepsilon^{(1)}) .
 \label{A2}
\end{eqnarray}

\bigskip
\noindent Notice that \eq{A1} implies $\tau_{i1}^{(1)}=0$ for
$i=1,\cdots,n$ in \eqref{2.7}, while \eq{A2} gives $\gamma=0$ in
\eqref{2.9}. It then follows that
\begin{eqnarray}
&T_{i1}=\sigma_{i1}T_{11}+\tau^{(0)}_{i1}, \quad
  i=1,\cdots,n, \label{2.10} \\[1mm]
&T:\varepsilon^{(1)}=\alpha T_{11}+\beta.
  \label{2.11}
\end{eqnarray}
This implies that the first column of $T$ is constant if and only if
$T_{11}$ in \eqref{2.5} is constant. Moreover, in this case,
$T:\varepsilon^{(1)}$ becomes a constant too. Therefore, under the
assumptions \eq{A1}  and \eq{A2}, for one-dimensional solutions of
the form \eqref{2.2}, we can reduce the original multi-dimensional
systems \eqref{I1} -- \eqref{I2} and \eqref{II1} -- \eqref{II2} both
with $b\equiv 0$ into the following one-dimensional systems:
\begin{eqnarray}
&T_{11}=D^{11}_{11}u^1_{x_1}-(D^{11}:\varepsilon^{(0)})
  -(D^{11}:\varepsilon^{(1)})v,
\label{n1d1}\\[1mm]
&v_t=(\alpha T_{11}+\beta-f'(v)+v_{x_1x_1})|v_{x_1}|,
\label{n1d2}
\end{eqnarray}
and
\begin{eqnarray}
&T_{11}=D^{11}_{11}u^1_{x_1}-(D^{11}:\varepsilon^{(0)})
  -(D^{11}:\varepsilon^{(1)})v,
\label{c1d1}\\[1mm]
&v_t=-\{(-f'(v)+v_{x_1x_1})_{x_1}|v_{x_1}|\}_{x_1},
\label{c1d2}
\end{eqnarray}
respectively, where $T_{11}$ is a constant.
These considerations are summarized as follows.

\begin{prop}\label{Prop2.2}
Suppose that the assumptions {\rm \eq{A1}} and {\rm \eq{A2}} are
met. Then the function $(v,u)$ given by \eq{2.2} is a solution to
system \eq{I1} -- \eq{I2} {\rm (}resp.   system \eq{II1} --
\eq{II2}{\rm )} with $b\equiv 0$ if and only if $(v,u^1)$ satisfies
the one-dimensional system \eq{n1d1} -- \eq{n1d2} {\rm (}resp. the
system \eq{c1d1} -- \eq{c1d2}{\rm )} with $T_{11}$ being a constant.
\end{prop}

\noindent{\bf Remark 1.}\ \ Since $T_{11}$ is regarded as a given
constant parameter,   both systems \eq{n1d1} -- \eq{n1d2} and
\eq{c1d1} -- \eq{c1d2} are completely decoupled. This means that we
solve \eq{n1d2} or \eq{c1d2} for $v$ independently, and then find
$u^1$ from \eq{n1d1} or \eq{c1d1}.


\section{Traveling wave solutions}
\setcounter{equation}{0}

In this section we are going to look for smooth traveling wave
solutions to the one-dimensional systems \eq{n1d1} -- \eq{n1d2} and
\eq{c1d1} --  \eq{c1d2} in the form
\begin{eqnarray}\label{3.1}
v = v(x_1-st), \qquad u^1_{x_1} = w(x_1-st),
\end{eqnarray}
where $s$ is a constant velocity to be determined. Here it is
assumed that $v(\xi)$ and $w(\xi)$ connect the constant states
$v_\pm$ ($v_-<v_+$) and $w_\pm$ (to be determined), respectively:
\begin{eqnarray}\label{infty}
v(\xi)\to v_\pm, \qquad  w(\xi)\to w_\pm
\end{eqnarray}
as $\xi\to\pm\infty$, and satisfies the requirement
\begin{eqnarray}
\label{positivity}
v_\xi>0
\end{eqnarray}
for $\xi\in\cR$. Here we have written $\xi=x_1-st$.

For the potential $f(v)$, we impose the following conditions:

\vskip0.2cm
 \noindent Assumption {\bf (B)}
\begin{eqnarray} && f(v)   \mbox{ is  a  smooth
double-well  potential  which  has  two  local  minima at}\ \non\\
&&v_-  \mbox{  and}\ v_+   \mbox{  with}\ v_-<v_+ \mbox{  and one
local maximum   at}\ v_* \mbox{  with}\ v_-<v_*<v_+ ,\quad \label{B}\\
&&  \mbox{  and  satisfies}\ f'(v)>0 \mbox{  for}\ v_-<v<v_* \mbox{
and} f'(v)<0  \mbox{  for}\
 v_*<v<v_+ .  \non
\end{eqnarray}

\bigskip
\noindent For simplicity, we assume that
\begin{equation}
\label{3.4}
\begin{split}
&f^{(k)}(v_+)=0\ \ {\rm for}\ \ 1\le k\le 2m_1-1, \quad
  f^{(2m_1)}(v_+)>0, \\[1mm]
&f^{(k)}(v_-)=0\ \ {\rm for}\ \ 1\le k\le 2m_2-1, \quad
  f^{(2m_2)}(v_-)>0,
\end{split}
\end{equation}
where $m_1$ and $m_2$ are positive integers.
A straightforward computation shows that for solutions of the form
\eqref{3.1},   system \eq{n1d1} --  \eq{n1d2} turns out to be
\begin{eqnarray}
&T_{11}=D^{11}_{11}w-(D^{11}:\varepsilon^{(0)})
  -(D^{11}:\varepsilon^{(1)})v,
  \label{I1c}   \\[1mm]
&-sv_\xi=(\alpha T_{11}+\beta-f'(v)+v_{\xi\xi})|v_{\xi}|,
  \label{I2c}
\end{eqnarray}
while   system \eq{c1d1}  -- \eq{c1d2} becomes
\begin{eqnarray}
&T_{11}=D^{11}_{11}w-(D^{11}:\varepsilon^{(0)})
  -(D^{11}:\varepsilon^{(1)})v,
  \label{II1c}    \\[1mm]
&-sv_\xi=-\{(-f'(v)+v_{\xi\xi})_{\xi}|v_{\xi}|\}_{\xi},
  \label{II2c}
\end{eqnarray}
where $T_{11}$ is assumed to be a given constant.

\vskip0.2cm {\noindent\bf Remark 2.} Since  systems \eq{n1d1} --
\eq{n1d2} and \eq{c1d1} --  \eq{c1d2}   are intended to describe
phase transitions,   it is necessary  to assume that the
nonlinearity $f(v)$ has at least two different minima, say $v_-, \
v_+$ and one maximum, say $v_*$.

We now consider the case that $v_- <   v_+$ (The case that $
v_->v_+$ can be treated in a similar manner, and we conclude that
traveling wave exists if and only if
 $v_\xi<0$ is satisfied). In this case,
traveling wave exists if and only if \eq{positivity} is true. In
this remark, we prove the necessity part. The sufficiency part is
the main task of the present article and will be proved after this
remark.

Assume  that there exists a smooth traveling wave $v$ and
\begin{eqnarray}
\mbox{ there exists at least one point, say }  \xi_0\in\cR , \mbox{
such
 that } v'(\xi_0 ) > 0 .
 \label{case1}
\end{eqnarray}
Otherwise, $v'(\xi)\le 0$ for all $\xi\in\cR$. It is easy to see
that such $v$ can not connect, $v_-$ and $v_+$, at minus and plus
infinity, respectively, provided that $v_-<v_+$.

Define
\begin{eqnarray}
{\cal N}_{\xi_0} & = & \{ N(\xi_0) \subset\cR\mid N(\xi_0)\mbox{ any
open simply connected set containing  } \xi_0,  v'(\xi ) \ge 0,\non\\
&&\quad\quad\quad\quad\quad\quad  \forall \xi\in N(\xi_0), \mbox{
and } v'(\xi )=0 \mbox{ only if } \xi \mbox{ is an inflection point}
\}
\end{eqnarray}
and a relation $\prec$: for any two elements $N_1,\, N_2\in {\cal
N}_{\xi_0}  $, $N_1 \prec N_2$ means $N_1\subset N_2$. By Zorn's
lemma (see, e.g. \cite{Yoshida}) we conclude that there exists a
maximal element $ N_{\rm max }(\xi_0)$ in ${\cal N}_{\xi_0}$, and
assume that there exist two numbers $\xi_-,\xi_+$  such that $
N_{\rm \max }(\xi_0)= (\xi_-,\xi_+)$. It is easy to see that
$v'(\xi_\pm)=0$, and we assume that
\begin{eqnarray}
v (\xi_-)= \underline{v},\quad v (\xi_+)= \bar v
  \label{vends}
\end{eqnarray}
 with   $\underline{v},\ \bar{v}\in\cR$  and $\underline{v} < \bar v  $. In $(\xi_-,\xi_+)$,
  equation \eq{I2c} can be reduced to
 an ordinary differential equation  of second order as follows
\begin{eqnarray}
  v_{\xi \xi}  =  f' (v) - (s+\alpha T_{11}+\beta) .
 \label{2ndpositive}
\end{eqnarray}

I) Suppose that both $\xi_+$ and $\xi_-$ are infinite.
Thus we can obtain $\bar v  =v_+$ and $\underline{v}=v_-$. Invoking
assumption  \eq{B}, we rewrite \eq{I2c} as
\begin{eqnarray}
 v_\xi  & = &   b(v) (v_+-v)^{m_1}(v-v_-)^{m_2} , \mbox{
 for } \xi\in  (\xi_-,\xi_+)
  \label{eq3a}
\end{eqnarray}
 here $b(v)$  is a smooth function. Thus from uniqueness theorem of ordinary
differential equations (see e.g. page 259, \cite{Smoller}),
 we exclude the case that $v_\xi(\xi_0)=0$ for some finite
 $\xi_0\in\cR$, and assumption \eq{positivity} follows.

II) Assume $ \bar v > v_+$ (or $\underline{v} < v_-$) at a   finite
point. Then there exists a point, say $\xi_* $, such that
$v(\xi_*)=v_+$. As in I) we know that $v$ satisfies \eq{eq3a} in a
small interval $[\xi_*,\xi^*)$ with $\xi^*>\xi_*$. This is
impossible again by uniqueness theorem of ordinary differential
equations.

III) Now we assume that either $\xi_+$ or $\xi_-$ is finite (We only
consider the case $\xi_+$  is finite, another case can be treated in
a similar way). Then there is a suitable number
$\xi_*>\xi_+$ such that \\
i) $v_\xi(\xi)<0$ for $\xi\in ( \xi_+, \xi_*)$;\\
ii)  $v (\xi)=\bar v$ for $\xi\in ( \xi_+, \xi_*)$,\\
as shown in Figure 1 (see page 16 of this article).


We first handle case ii). Since we look for a smooth solution $v$,
the left  limits at $\xi_+$ of
 derivatives  up to second order (which is enough here) of $v$  must be equal
 to the right ones, respectively. It is easy to
 compute that
 $v_{\xi \xi}(\xi_++0)=0$. Thus taking
 limit for both sides of \eq{2ndpositive} yields that
\begin{eqnarray}
  0= v_{\xi \xi} (\xi_--0) =  f' (\bar v) - (s+\alpha T_{11}+\beta) .
 \label{2ndpositiveLeft}
\end{eqnarray}
So this means $v\equiv\bar v$ is a solution to \eq{2ndpositive}.
This violates uniqueness theorem of ordinary differential equations
in $(\xi_-,\xi_+]$. Therefore case ii) could  not happen.

It remains to consider case i). In $( \xi_+, \xi^*)$, equation
\eq{I2c} is equivalent to
\begin{eqnarray}
  v_{\xi \xi}  =  f' (v) - (-s +\alpha T_{11}+\beta) .
 \label{2ndnegative}
\end{eqnarray}
 Thus from \eq{2ndnegative} and \eq{2ndpositive} one
has
$$
f' (\bar v) - (-s +\alpha T_{11}+\beta) = f' (\bar v) - ( s +\alpha T_{11}+\beta)
$$
which implies that
\begin{eqnarray}
  s  =  0 ,
 \label{s}
\end{eqnarray}
hence if $s\not=0$, then case i) can not be true. Suppose that
\eq{s}  is satisfied, then \eq{2ndnegative} and \eq{2ndpositive}
become
\begin{eqnarray}
  \frac12(v_\xi)^2  &=& \Phi(v;c):= f (v) - ( \alpha T_{11}+\beta)v + c,  \mbox{
  for }  \xi\in ( \xi_-, \xi_+),  \label{eqmaxP} \\
 \frac12(v_\xi)^2   &=& \Phi(v;c_1):=  f (v) - ( \alpha T_{11}+\beta ) v + c_1 , \mbox{
 for }  \xi\in  ( \xi_+, \xi^*).
 \label{eqmaxN}
\end{eqnarray}
We can not see immediately that equations \eq{eqmaxP} and
\eq{eqmaxN} can be rewritten in the form of \eq{eq3a} that meets the
conditions of uniqueness theorem of first order ordinary
differential equations. Note that the right hand side of \eq{eqmaxP}
differs, by a constant, from that of \eq{eqmaxN}, namely $c,\ c_1$
may be different. But at $\xi=\xi_+$, there hold $v_\xi=0,\ v=\bar
v$. Hence, it follows from \eq{eqmaxP} and \eq{eqmaxN} that
$\Phi(\bar v;c) =0 = \Phi(\bar v;c_1)$, so $c=c_1$. That is to say
$v$ satisfies the same equation in both $( \xi_+, \xi^*)$ and $(
\xi_-, \xi_+)$. By extension, it is also  true that $v$ satisfies
the same equation in any interval contained in $( -\infty, \infty)$.
Thus $c=c_1=-(f (v_\pm) - (s+\alpha T_{11}+\beta)v_\pm)$, and
\eq{I2c} can be rewritten as \eq{eq3a}. We can apply again
uniqueness theorem of ordinary differential equations and conclude
it is impossible that $v_\xi(\xi_0)=0$ for any finite $\xi_0$. This
leads to a contradiction since we assumed that  $v_\xi(\xi_+)=0$ and
$\xi_+$ is finite. So we again obtain that $v_\xi$ can not change
its sign even \eq{s} holds, and must be positive since we assume
$v_-<v_+$.

\bigskip
We now turn back to the investigation of traveling waves. The
remaining of this section will be divided into two parts
corresponding to systems \eqref{I1c}  --  \eqref{I2c} and
\eqref{II1c}  --  \eqref{II2c}, respectively.

\bigskip
\noindent {\it Part 1.}\ We first study the non-conservative system
\eq{I1c}  --  \eq{I2c}. Letting $\xi\to\pm\infty$ in \eq{I1c} and
recalling that $T_{11}$ is constant, we obtain
\begin{equation*}
T_{11}=D^{11}_{11}w_\pm-(D^{11}:\varepsilon^{(0)})
  -(D^{11}:\varepsilon^{(1)})v_\pm,
\end{equation*}
which gives
\begin{eqnarray}
D^{11}_{11}(w_+-w_-)=(D^{11}:\varepsilon^{(1)})(v_+-v_-).
\label{jumpW}
\end{eqnarray}
On the other hand, since we have assumed \eq{positivity}, one can
obtain from \eq{I2c} that
\begin{equation*}
\begin{split}
-sv_\xi
&=(\alpha T_{11}+\beta-f'(v)+v_{\xi\xi})v_{\xi} \\[1mm]
&=\big\{(\alpha T_{11}+\beta)v-f(v)+\frac12
(v_{\xi })^2\big\}_{\xi}.
\end{split}
\end{equation*}
Integrating it with respect to $\xi$ yields
\begin{eqnarray}\label{ode1}
\frac12(v_{\xi })^2
=f(v)-(s+\alpha T_{11}+\beta)v-A=:g(v),
\end{eqnarray}
where $A$ is a constant.
Letting $\xi\to\pm\infty$ and using \eq{infty} gives
\begin{eqnarray*}
f(v_\pm)-(s+\alpha T_{11}+\beta)v_\pm-A=0,
\end{eqnarray*}
from which follows that
\begin{eqnarray}
s+\alpha T_{11}+\beta=\frac{f(v_+)-f(v_-)}{v_+-v_-}.
 \label{speed}
\end{eqnarray}
As a consequence, we see that $g(v)$ defined in \eqref{ode1}
can be rewritten as
\begin{eqnarray}
g(v)=f(v)-f(v_\pm)-\frac{f(v_+)-f(v_-)}{v_+-v_-}(v-v_\pm).
\label{gv1}
\end{eqnarray}
Since $g(v)\ge 0$ by \eq{ode1}, we get
\begin{eqnarray*}
f(v)-f(v_\pm)\ge\frac{f(v_+)-f(v_-)}{v_+-v_-}(v-v_\pm)
\end{eqnarray*}
for all $v$ such that $v_-<v<v_+$. We divide this inequality by
$v-v_\pm$, where $v-v_+<0<v-v_-$, and take the limits as
$\xi\to\pm\infty$. Since $v(\xi)\to v_\pm$ as $\xi\to\pm\infty$,
we deduce that
\begin{eqnarray*}
f'(v_+)\le\frac{f(v_+)-f(v_-)}{v_+-v_-}\le f'(v_-).
\end{eqnarray*}
By the assumption \eq{B} or \eq{3.4} we have $f'(v_\pm)=0$, and thus
we infer that
\begin{eqnarray}
f(v_+)=f(v_-),
\label{shock1a0}
\end{eqnarray}
which together with \eq{speed} gives
\begin{eqnarray}
s=-(\alpha T_{11}+\beta).
\label{speedfinal}
\end{eqnarray}
Thus the velocity $s$ of the traveling wave solution is determined.

Consequently, our equation \eq{ode1} can be simplified to
\begin{eqnarray}\label{3.15}
\frac12(v_{\xi })^2=f(v)-f(v_\pm).
 \label{reduced}
\end{eqnarray}
By \eq{shock1a0} and the assumptions \eq{B} and \eq{3.4}, we can
write
\begin{eqnarray}
f(v)-f(v_\pm)=\frac12a(v)(v_+-v)^{2m_1}(v-v_-)^{2m_2}
 \label{rhs}
\end{eqnarray}
with $a(v)>0$ being a smooth function.
Therefore \eq{reduced} can be rewritten as
\begin{equation}
\label{320}
v_{\xi}=a(v)^{1/2}(v_+-v)^{m_1}(v-v_-)^{m_2}.
\end{equation}
This equation admits a smooth solution $v(\xi)$ satisfying
\eq{infty} and \eq{positivity}. Moreover, the solution verifies
the following decay estimates for $\xi\to\pm\infty$:
\begin{eqnarray}
\big|\partial^k_\xi(v(\xi)-v_+)\big|
\le\left\{
 \begin{array}{lll}
 Ce^{-c|\xi|} & {\rm if }\ \ m_1=1, \\[1mm]
 C(1+|\xi|)^{-\nu_1-k} & {\rm if }\ \ m_1>1
 \end{array}
\right.
\label{est+}
\end{eqnarray}
for $\xi\to\infty$, and
\begin{eqnarray}
\big|\partial^k_\xi(v(\xi)-v_-)\big|
\le\left\{
 \begin{array}{lll}
 Ce^{-c|\xi|} & {\rm if }\ \ m_2=1, \\[1mm]
 C(1+|\xi|)^{-\nu_2-k} & {\rm if }\ \ m_2>1
 \end{array}
\right.
\label{est-}
\end{eqnarray}
for $\xi\to -\infty$, where $k$ is a non-negative integer,
$\nu_1=1/(m_1-1)$, $\nu_2=1/(m_2-1)$, and $C$ and $c$ are positive
constants.

Thus we have proved the following theorem:

\begin{theorem}
{\rm\bf (Non-conservative case)}\ Suppose that the assumptions {\rm
\eq{B}} and \eq{3.4} are satisfied. Then the one-dimensional system
\eq{n1d1} --  \eq{n1d2} admits a smooth traveling wave solution
$(v,u^1_{x_1})=(v,w)(x_1-st)$ satisfying \eq{infty} and
\eq{positivity} if and only if the constants $s$, $v_\pm$ and
$w_\pm$ satisfy the relations \eq{speedfinal}, \eq{shock1a0} and
\eq{jumpW}. In particular, two local minima $f(v_\pm)$ must coincide
in this case. Moreover, such a traveling wave solution is unique up
to a translation in $\xi$ and verifies the decay estimates \eq{est+}
and \eq{est-} for $\xi\to\pm\infty$, respectively.
\end{theorem}

\medskip
\noindent {\it Part 2.}\ Next, we discuss the conservative system
\eq{II1c} --  \eq{II2c}. In this case, we still have \eq{jumpW}.
Also, integrating \eq{II2c} gives
\begin{eqnarray}
-sv=-(-f'(v)+v_{\xi\xi})_\xi|v_\xi|+A_1,
 \label{Ode1}
\end{eqnarray}
where $A_1$ is a constant. Letting $\xi\to\pm\infty$,
we then obtain $A_1=-sv_\pm$, which implies
\begin{eqnarray}
s=0.
\label{zeroSpeed}
\end{eqnarray}
Whence we have $A_1=0$ so that \eq{Ode1} becomes
$(-f'(v)+v_{\xi\xi})_\xi|v_\xi|=0$.
Since $v_\xi>0$ by \eq{positivity}, we get
$(-f'(v)+v_{\xi\xi})_\xi=0$ and hence
\begin{eqnarray}
v_{\xi\xi }-f'(v)=A_2,
\label{Ode2}
\end{eqnarray}
where $A_2$ is a constant. Letting again $\xi\to\pm\infty$, we find
that $A_2=-f'(v_\pm)=0$, where we have used the assumption \eq{B} or
\eq{3.4}. Consequently, \eq{Ode2} becomes
\begin{eqnarray}
v_{\xi\xi }-f'(v)=0.
\label{Ode2a}
\end{eqnarray}
Multiplying \eq{Ode2a} by $v_\xi$, we get
$\big\{\frac12(v_{\xi})^2-f(v)\big\}_\xi=0$.
We integrate this equation in $\xi$ to obtain
\begin{eqnarray}
\frac12(v_{\xi})^2=f(v)-A_3,
 \label{Ode2c}
\end{eqnarray}
where  $A_3$ is a constant. Letting $\xi\to\pm\infty$, we see that
$A_3=f(v_\pm)$ and hence the relation \eq{shock1a0} must be
satisfied. Consequently, we find that \eq{Ode2c} is just the same as
\eq{3.15} in Part 1. Therefore the same arguments as in Part 1 prove
the following theorem.

\begin{theorem}
{\rm\bf (Conservative case)}\ Suppose that the assumptions {\rm
\eq{B}} and \eq{3.4} are satisfied. Then the one-dimensional system
\eq{c1d1} --  \eq{c1d2} admits a smooth traveling wave solution
$(v,u^1_{x_1})=(v,w)(x_1-st)$ satisfying \eq{infty} and
\eq{positivity} if and only if the constants $s$, $v_\pm$ and
$w_\pm$ satisfy the relations \eq{zeroSpeed}, \eq{shock1a0} and
\eq{jumpW}. In particular, two local minima $f(v_\pm)$ must coincide
also in this case. In particular, necessarily, $s=0$, and the
traveling wave   is actually a standing wave solution. Moreover,
such a standing wave solution is unique up to a translation in $\xi$
and verifies the decay estimates \eq{est+} and \eq{est-} for
$\xi\to\pm\infty$, respectively.
\end{theorem}


\section{Comparison with the Allen-Cahn/Cahn-Hilliard equations
coupled with linear elasticity}
\setcounter{equation}{0}

In this section we are going to compare our results with those for
the Allen-Cahn  and the Cahn-Hilliard equations,   coupled with the
linear elasticity system \eq{I1} with constitutive relation
\eq{CR1}. We shall consider only
 quasi-static models.  The Allen-Cahn and   Cahn-Hilliard equations
coupled with \eq{I1} are, respectively, as follows
\begin{eqnarray}
&-\,{\rm div}_x\,T=b,
  \label{I1ac}   \\[1mm]
&v_t=T:\varepsilon^{(1)}-f'(v)+\Delta_xv,
 \label{I2ac}
\end{eqnarray}
and
\begin{eqnarray}
&-\,{\rm div}_x\,T=b,
  \label{II1ch}    \\[1mm]
&v_t=-{\rm div}_x\big\{\nabla_x
  (T:\varepsilon^{(1)}-f'(v)+\Delta_xv)\big\}.
  \label{II2ch}
\end{eqnarray}
Here we made the same assumptions \eq{CR1} and \eq{CR2} for the
systems.

As we will see, there are  interesting differences  between the
result of existence of traveling waves for system \eq{I1ac} --
\eq{I2ac} and those for system \eq{I1a} -- \eq{I2a}.

The well-posedness of systems \eq{I1ac} -- \eq{I2ac} and \eq{II1ch}
-- \eq{II2ch}, has been studied, in particular,  in
\cite{BW05,Garcke03}.

Their is also an extensive literature  on traveling waves to the
single Allen-Cahn equation \eq{I2ac} with various nonlinearities: we
refer the reader to the articles, e.g., Aronson and Weinberger
\cite{AW78}, Berestycki and Hamel \cite{BH}, \cite{FM77}, Ninomiya
and Taniguchi \cite{NT,Taniguchi}, Taniguchi \cite{Taniguchi07}, or
the book by Fife \cite{Fife79}, etc.
  However, there are very few results, to the knowledge of the authors, for the
existence of traveling waves for the complete systems \eq{I1ac} --
\eq{I2ac} and \eq{II1ch} -- \eq{II2ch}.

\subsection{Traveling waves for the Allen-Cahn/Cahn-Hilliard models}

We divide this subsection into two parts.

\bigskip
 \noindent {\it Part 1.}\ System \eq{I1ac} --
\eq{I2ac}.  We now turn to consider the non-conservative case. We
assume that
$$
b=0,\qquad v=v(x_1,t), \qquad u=(u^1(x_1,t),0,\cdots,0),
$$
and that the same assumptions \eq{A1} -- \eq{A2}  are true to system
\eq{I1ac} -- \eq{I2ac} as in Section~2. Then  similar computations
yield that   \eq{I1ac} -- \eq{I2ac} can be reduced to
\begin{eqnarray}
v_t=v_{x_1x_1}+\mu-f'(v),
\label{single}
\end{eqnarray}
where $\mu=\alpha T_{11}+\beta$ is a   constant.

To prove the existence of traveling wave solutions to  \eq{single},
we recall one of the existence theorems (see for instance,
\cite{Kanel,FM77} and also the references cited therein) for the
following single Allen-Cahn equation
\begin{equation}\label{4.6}
v_t=v_{x_1x_1}+g(v).
\end{equation}

Notice that eq. \eq{single} can be written in the form of \eq{4.6}
by defining
$$g(v)=\mu-f'(v) .
$$
 For the convenience of the readers, we state the following theorem
for equation \eq{4.6}.

\begin{theorem}
{\rm (Kanel \cite{Kanel}, and Theorem 3.1 in \cite{FM77})}\
Let $v_-<v_*<v_+$ and let $g\in C^1(\cR)$ satisfy $g(v_\pm)=0$,
$g'(v_\pm)<0$, $g(u)<0$ for $v_-<v<v_*$, and $g(v)>0$ for
$v_*<v<v_+$. Then there exists a unique (up to a translation)
monotone traveling wave solution to \eq{4.6}, which connects the end
states $v_-$ and $v_+$.
\end{theorem}

In the following, when we speak of a ``traveling wave over
$[\alpha,\beta]$ with velocity $s_{\alpha,\beta}$", we shall mean a
solution of
 equation \eq{4.6} with the given $s_{\alpha,\beta}$, which is
 positive in $(\alpha,\beta)$  and vanishes at
$ \alpha$ and $\beta $. We now cite the following theorem
\begin{theorem}  {\rm (Theorem 2.7 in \cite{FM77})}  Assume that
 $g\in C^1(\cR)$  with $g(v_\pm) = 0$,
%
and let there exist a traveling wave over $[v_-,v_* ]$ with velocity
$s_{v_-,v_*}$, and one over $[v_*, v_+]$ with velocity $s_{v_*,v_+}
<s_{v_-,v_*}$. Then there exists a traveling wave over $[v_-, v_+]$
with velocity $s_{v_-,v_+}$ satisfying
$$
s_{v_*,v_+} <s_{v_-,v_+} <s_{v_-,v_*}.
$$
\end{theorem}

\noindent {\bf Remark 3.}\ In the original version of this theorem,
$v_-$ and $v_+$ are assigned special values which are $v_-=0$ and
$v_+=1$.

\bigskip
\noindent {\bf Remark 4.}\  It has been pointed out in
\cite{Taniguchi07},  p. 320,  that, based upon Theorem~4.2, one can
  find a function $g$ such  that for such $g$,  eq.
\eq{4.6} has a traveling wave over $[v_-,v_* ]$ with velocity
$s_{v_-,v_*}$, and one over $[v_*, v_+]$ with velocity $s_{v_*,v_+}
\ge s_{v_-,v_*}$.  The condition   $s_{v_*,v_+} <s_{v_-,v_*}$
 required by Theorem~4.2 is violated,
  so there exists {\it no} traveling wave that connects two states $v_\pm$,
   to eq. \eq{4.6}.

\bigskip
From Theorem~4.1
we know that to guarantee the existence
of traveling waves, the following necessary condition must be
satisfied
\begin{eqnarray}
g(v_\pm)=0.
 \label{FM0}
\end{eqnarray}
Invoking that $v_\pm$ are two local minima of $f(v)$ by assumption
\eq{B}, we have $f'(v_\pm)=0$, and thus there must hold
$$
0=g(v_\pm)=\mu-f'(v_\pm)=\mu
$$
hence,
\begin{eqnarray}
\mu=\alpha T_{11}+\beta=0,
 \label{FM}
\end{eqnarray}
which means that there is no contribution from the elastic energy to
the equation \eq{single}.

Suppose now that $\mu=0$ is satisfied.
It is easy to see that  \eq{single}  is reduced to   \eq{4.6} with
$g(v)=-f'(v)$. Equation \eq{4.6}  is just the single Allen-Cahn
equation and many existence theorems are applicable to it. For
example, under the assumption \eq{B}, it is easy to check that all
the conditions in Theorem~4.1 are met for $g(v)=-f'(v)$. Then
applying Theorem~4.1, we
  conclude that there exists a traveling wave solution to
\eq{single} with $\mu=\alpha T_{11}+\beta=0$, and this traveling
wave connects two states $v_-$ and $v_+$. But there exist some
nonlinear functions $f(v)$ for which there exists no traveling waves
to the single Allen-Cahn equation \eq{4.6} with $g(v)=-f'(v)$, as
pointed out in Remark~3.

\bigskip
\noindent {\it Part 2.}\ System \eq{II1ch} -- \eq{II2ch}. Next we
consider the conservative model \eq{II1ch} -- \eq{II2ch}. Assuming
that
$$
b=0,\qquad v=v(x_1,t), \qquad u=(u^1(x_1,t),0,\cdots,0),
$$
and that the same assumptions \eq{A1} - \eq{A2}  are true, we can
reduce system \eq{II1ch} -- \eq{II2ch} to
\begin{equation}\label{single1}
v_t=-(v_{x_1x_1}-f'(v))_{x_1x_1}.
\end{equation}

Note that the gradient term in \eq{Ode1} for the modified
Cahn-Hilliard model \eq{II1a} -- \eq{II2a} does not influence
essentially the results after  equation \eq{Ode2a} in Section~3.
Therefore, same arguments can be carried out for eq. \eq{single1}
and thus we conclude easily that there is no traveling wave solution
with non-zero speed for the Cahn-Hilliard model \eq{II1ch} --
\eq{II2ch} too.


\subsection{Comparison}

From Theorems~3.2, 3.3 and the arguments in Subsection~4.1, we
   are now able to draw the following conclusions.

\medskip
\noindent {\bf Conclusions:}\

\noindent{\bf A)} Non-conservative case:

\noindent \quad (1) Suppose that $\alpha T_{11}+\beta\neq 0$ is met,
which means that the elastic energy contributes to the total free
energy;
Then the elastic effect prevents the formation of traveling waves to
the Allen-Cahn model \eq{I1ac} -- \eq{I2ac}. In contrast, the
gradient term overcomes the effect due to the elastic energy and so
there still exists a traveling wave to the modified Allen-Cahn model
\eq{I1a} -- \eq{I2a}.

\noindent \quad (2) Assume that there holds $\alpha T_{11}+\beta=0$.
In this case there is no elastic energy. Under suitable assumptions
on the nonlinearity, we can prove that there exists a traveling wave
solution with non-zero speed both for the Allen-Cahn model \eq{I1ac}
-- \eq{I2ac}
and also for the modified Allen-Cahn one \eq{I1a} -- \eq{I2a}.

However, in some other cases (e.g. for a function $g$ such that the
condition  $s_{\alpha,v_+} <s_{v_-,\alpha}$ in Theorem~4.2 is not
met, see Remark~3), there exists no traveling wave solutions for
models \eq{I1ac} -- \eq{I2ac} and \eq{I1a} -- \eq{I2a}.

\medskip
\noindent{\bf B)} Conservative case:

\noindent There is no essential difference from the viewpoint of the
existence of traveling waves since both   the Cahn-Hilliard model
\eq{II1ch} -- \eq{II2ch} and the modified one \eq{II1a} -- \eq{II2a}
admit only a standing wave.

\bigskip
\noindent{\bf Acknowledgement.}\ The authors would like to express
their sincere thanks to Prof. H.-D.~Alber and Prof. E. Zuazua for
 helpful discussions. This work has been partially supported by Grant
MTM2008-03541 of the MICINN (Spain).


\includegraphics[scale=0.75]{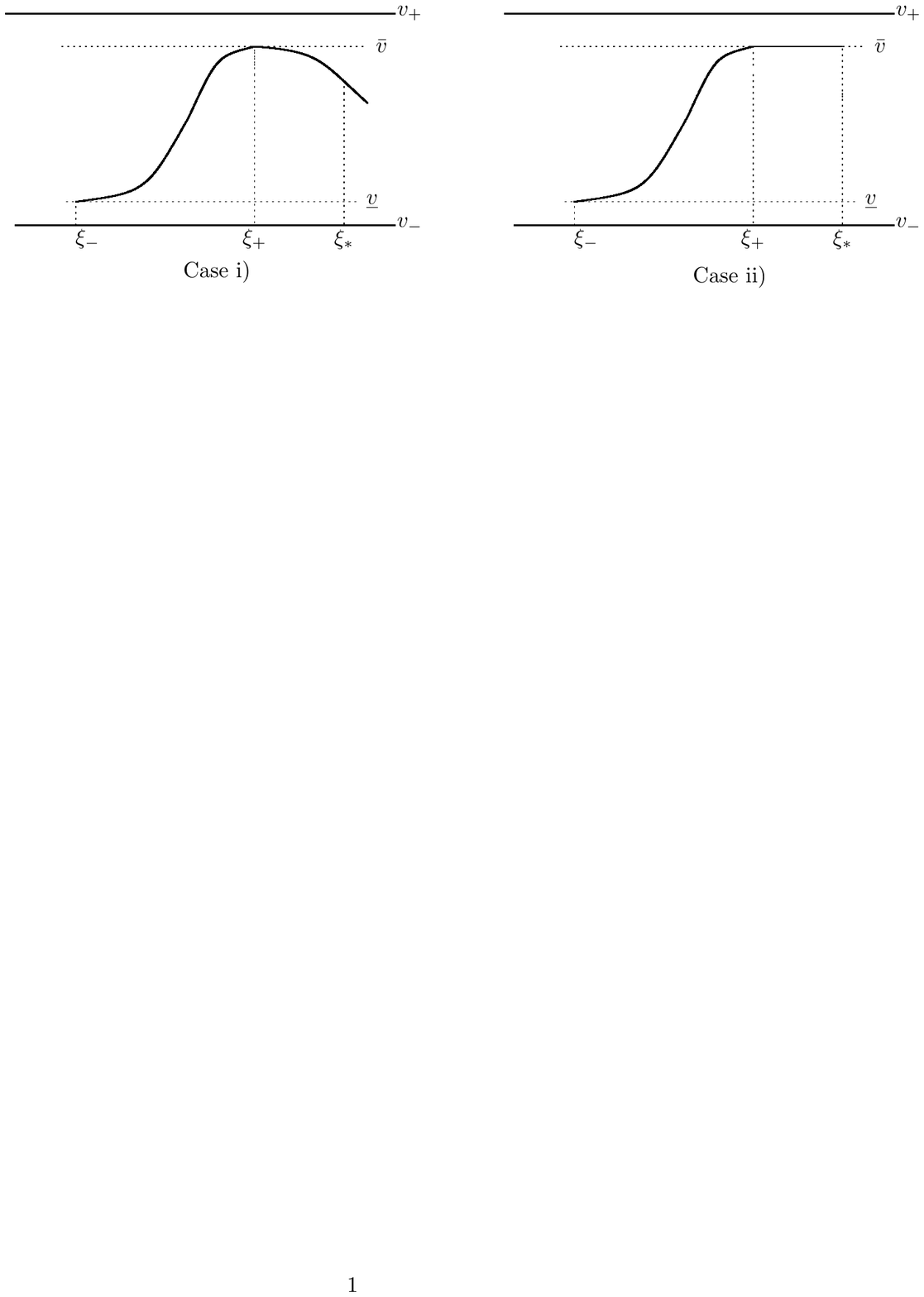}

\vskip-15cm\hskip8.5cm{Figure 1.}


\begin{thebibliography}{99}
\itemsep0pt

\bibitem{AK90} Abeyaratne, R. and Knowles, J. (1990) On the driving
traction acting on a surface of strain discontinuity in a continuum,
 {\it J. Mech. Phys. Solids} {\bf 38},  345--360.

\bibitem{Alber00} Alber, H.-D. (2000)  {Evolving microstructure and
 homogenization,} {\it Continuum. Mech. Thermodyn.} {\bf 12},
 235--287.

\bibitem{Alber04} Alber, H.-D. and Zhu, P. (2007) {Evolution of phase
 boundaries by configurational forces,}    {\it Archive
 Rat. Mech. Anal.}   {\bf 185}, 235--286.

\bibitem{Alber06} Alber, H.-D., and Zhu, P. (2006) {Solutions to a
 model with nonuniformly parabolic terms for phase evolution driven by
 configurational forces,} {\it SIAM J. Appl. Math.} {\bf 66}, No. 2,
 680--699.

\bibitem{Alber08} Alber, H.-D., and Zhu, P. (2008) {  Solutions to a
 model   for interface motion  by interface diffusion,} {\it Proc. Roy. Soc. Edinburgh} {\bf 138A},
 923--955.

\bibitem{Alber08a} Alber, H.-D., and Zhu, P. (2008a) {  Interface motion  by interface diffusion driven
by bulk energy: justification of a diffusive interface model,}
submitted to {\it Conti. Mech. and Thermodyn.}.

\bibitem{Allen79}  Allen, S. and  Cahn, J.   (1979)   A microscopic theory for antiphase boundary motion
and its application to antiphase domain coarsening, {\it Acta Met.}
{\bf 27},  1084–1095.

\bibitem{AW78} Aronson, D. and Weiberger, H. (1978) Multidimensional
nonlinear diffusions arising in population genetics, {\it Adv.
Math.} {\bf 30},   33--76.

\bibitem{BH} Berestycki, H. and Hamel, F. (2007) Generalized travelling
waves for reaction-diffusion equations, {\it Perspectives in
Nonlinear Partial Differential Equations}. In honor of H. Brezis,
{\it Contemp. Mathematics}, {\bf 446}, Amer. Math. Soc., Providence,
RI.

\bibitem{BW05} Blesgen, T. and Weikard, U. (2005) Multi-component Allen-Cahn equation for
elastically stressed solids, {\it Electronic J. Diff. Eq.s,} {\bf
2005}, No. 89,  1--17.

\bibitem{Cahn59} Cahn, J.   (1959)  Free Energy of a Nonuniform
System. II. Thermodynamic Basis, {\it J. Chem. Phys.}, {\bf   30},
1121--1124.

\bibitem{Cahn59a}   Cahn, J. and   Hilliard, J. (1959)  Free Energy of a Nonuniform
System. I. Interfacial Free Energy, {\it J. Chem. Phys.}, {\bf  28},
258--267; Free Energy of a Nonuniform System. III. Nucleation in a
Two-Component Incompressible Fluid,  {\it J. Chem. Phys.}, {\bf 31},
688-699.


\bibitem{CT94}  Cahn, J.,  and   Taylor, J.   (1994) Surface motion by surface
 diffusion.  {\it Acta Metall. Mater.}  {\bf 42}, No. 4,
 1045--1063.

\bibitem{Fife79} Fife, P. (1979) {\it Mathematical aspects of reacting and diffusing
systems}, Springer Verlag.

\bibitem{FM77} Fife, P. and McLeod, J.   (1977) The approach of
solutions of nonlinear diffusion equations to traveling front
solutions, {\it Arch. Rati. Mech. Anal.}, {\bf 65}, 335--361.

\bibitem{Gurtin00} Gurtin, M. (2000) {\it Configurational Forces as Basic Concepts of Continuum
Physics}, Springer Verlag, New York.

\bibitem{Garcke03} Garcke, H. (2003)   On the Cahn-Hilliard system
with elasticity, {\it Proc. Roy. Soc. Edinburgh},  {\bf 133 A},
307--331.

\bibitem{Hornbogen} Hornbogen, E. and Warlimont, H. (2001) {\it Metallkunde}, 4th ed.,
Springer-Verlag.

\bibitem{Kanel} {Kanel}, Y. (1962) On the stabilization of solutions
of the Cauchy problem for equations arising in the theory of
combustion, {\it  Mat. Sbornik} {\bf  59}, 245--288.

\bibitem{Mullins}   Mullins, W.  (1957) Theory of thermal grooving. {\em
 J. Appl. Phys.} {\bf 28}, No. 3, 333--339.

\bibitem{NT} Ninomiya, H. and Taniguchi, M. (2005) Existence and
global stability of traveling curved fronts in the Allen-Cahn
equations, {\it J. Diff. Eq.} {\bf 213},
  204--233.

\bibitem{Taniguchi} Ninomiya, H. and Taniguchi, M. (2006) Global stability of traveling
curved fronts in the Allen-Cahn equations, {\it  Disc. Conti. Dyna.
Syst.} {\bf 15}, No. 3,   819--832.

\bibitem{Pego} Pego, R. (1989) {Front migration in the nonlinear
  Cahn-Hilliard equation,} {\it Proc. R. Soc. Lond.}  {\bf 422A},
  261--278.

\bibitem{Smoller} Smoller, J. (1983) {\it Shock waves and
reaction-diffusion equations}, Springer Verlag, New York

\bibitem{Taniguchi07}   Taniguchi, M. (2007) Traveling fronts of pyramidal shapes in the
Allen-Cahn equations, {\it  SIAM J. Math. Anal.} {\bf 39}, No. 1,
319--344.

\bibitem{TC94} Taylor, J.  and  Cahn, J.   (1994) Linking anisotropic sharp
 and diffuse surface motion laws via gradient flows. {\it
 J. Stat. Phys.} {\bf 77}, Nos. 1/2, 183--197.

\bibitem{Yoshida} Yoshida, K. (1980) {\it Functional analysis}  sixth edition,
 Springer Verlag, Berlin Heidelberg New York


\end{thebibliography}
\end{document}